\documentclass[reqno,10pt]{amsart}
\usepackage{amscd,amssymb}
\usepackage{latexsym}
\usepackage[all]{xy}

\def\into{\hookrightarrow}

\def\rra{\rightrightarrows}

\def\longto{\dashrightarrow}

\def\toisom{\widetilde{\to}}

\def\.{,\dots ,}
\def\wt{\widetilde}
\def\wh{\widehat}
\def\ol{\overline}

\def\Var{{\rm Var}}

\def\Spf{{\rm Spf}}

\def\Spec{{\rm Spec}}
\def\Proj{{\rm Proj}}

\def\Bl{{\rm Bl}}

\def\hatBl{{\rm \wh{Bl}}}
\def\red{{\rm red}}

\def\reg{{\rm reg}}

\def\sing{{\rm sing}}

\def\Id{{\rm Id}}

\def\Hom{{\rm Hom}}
\def\Der{{\rm Der}}

\def\calHom{{\mathcal Hom}}

\def\small{{\rm small}}

\def\sh{{\rm sh}}

\def\cha{{\rm char}}

\def\bfA{{\bf A}}

\def\bfC{{\bf C}}

\def\bfN{{\bf N}}

\def\bfQ{{\bf Q}}

\def\bfZ{{\bf Z}}

\def\gtC{{\mathfrak C}}

\def\gtI{{\mathfrak I}}

\def\gtV{{\mathfrak V}}

\def\gtX{{\mathfrak X}}

\def\gtZ{{\mathfrak Z}}

\def\calC{{\mathcal C}}

\def\calE{{\mathcal E}}
\def\calF{{\mathcal F}}
\def\calG{{\mathcal G}}

\def\calI{{\mathcal I}}

\def\calM{{\mathcal M}}
\def\calN{{\mathcal N}}
\def\calO{{\mathcal O}}
\def\calP{{\mathcal P}}

\def\oX{{\ol X}}
\def\oY{{\ol Y}}

\def\on{{\ol n}}

\def\tilX{{\wt X}}
\def\tilY{{\wt Y}}

\def\hatA{{\wh A}}

\def\hatK{{\wh K}}

\def\hatR{{\wh R}}

\def\hatX{{\wh X}}

\def\hatcalF{{\wh\calF}}

\def\hatcalI{{\wh\calI}}

\def\hatcalO{{\wh\calO}}

\def\R+*{{\bf R^*_+}}

\newtheorem{theor}{Theorem}[subsection]

\newtheorem{lem}[theor]{Lemma}

\theoremstyle{definition}

\newtheorem{rem}[theor]{Remark}
\newtheorem{exam}[theor]{Example}

\begin{document}

\title{Absolute desingularization in characteristic zero}
\author{Michael Temkin}
\address{\tiny{Department of Mathematics, Institute for Advanced Study, Princeton, NJ 08540, USA}}
\email{\scriptsize{temkin@ias.edu}}
\thanks{I am grateful to E. Bierstone and P. Milman for useful discussions and correcting an example in a
preliminary version of the paper. I would like to thank the ICMS and the organizers of the workshop on Motivic
Integration and its Interactions with Model Theory and Non-Archimedean Geometry held in May 2008 for the
hospitality and excellent working conditions during the workshop. This work was finished during my visit at the
Institute for Advanced Study at Princeton when I was supported by the NFS grant DMS-0635607.}

\maketitle

\section{Introduction}

\subsection{Preamble}
This paper is an expository lecture notes originally based on a lecture on the results of \cite{temdes} given by
the author at the workshop on Motivic Integration in May 2008, at ICMS, Edinburgh. Since a substantial progress
was done since May 2008, it seemed natural to include the new results of \cite{bmt}, \cite{temnon} and
\cite{tememb} in this exposition. We will mainly concentrate on the functorial non-embedded desingularization
constructed in \cite{temnon} because it seems that the results of \cite{tememb} on the embedded case can be
improved further. We pursue expository goals, so we will concentrate on explaining the results and the main
ideas of our method and we will refer to the cited papers for proofs and technical details. Also, we try to
include more examples and general remarks than in a pure research paper. Thus, this survey can serve as a
companion to or a light version of \cite{temdes} and \cite{temnon}. I would like to warn the reader that the
current situation described in the paper can change soon (similarly to the change since 2008), but this is
always a danger with a survey on an active research area.

\subsection{The history}
In 1964 Hironaka proved many fundamental desingularization results including strong desingularization of
algebraic varieties in characteristic zero. The latter means that any reduced variety of characteristic zero can
be modified to a smooth one by successive blow ups along nowhere dense smooth centers. Hironaka's method was
extremely difficult for understanding (due to a complicated inductive structure of the proof), and perhaps the
main reason for this was that his method was not constructive, canonical or functorial, unlike many new proofs.
In particular, unlike the new proofs, Hironaka could not work within the category of varieties since some
arguments with formal completions were involved. Probably for this reason, Hironaka proved his desingularization
for all schemes of finite type over local rings $R$ with regular completion homomorphism $R\to\hatR$.

A year later, Grothendieck introduced quasi-excellent (or {\em qe}) schemes in \cite[$\rm IV_2$, \S7.9]{ega} in
order to provide the most general framework for desingularization. Grothendieck observed that the schemes
studied by Hironaka were schemes of finite type over a local qe scheme $k$, and proved that if any integral
scheme of finite type over a base scheme $k$ admits a desingularization in the weakest possible sense then $k$
is qe. Grothendieck conjectured that the converse is probably true (i.e. any integral qe scheme admits a
desingularization) and claimed without proof that the conjecture holds true for noetherian qe schemes over
$\bfQ$ as can be proved by Hironaka's method. The latter claim was never checked in published literature, and
its status is unclear until now. Nevertheless, this fact was occasionally used by other mathematicians, for
example for desingularizing affinoid spaces over $\bfQ_p$.

In \cite{temdes} the author proved that indeed, any noetherian integral qe scheme over $\bfQ$ admits a
desingularization. Moreover, the regular locus of the scheme can be preserved by the desingularization and one
can also resolve a closed subscheme to a normal crossings divisor. The construction of \cite{temdes} uses any
desingularization of varieties as a black box input, but modifies it significantly. In particular, many good
properties were lost in the resulting algorithm, including functoriality and regularity of the centers. Very
recently the method was strengthened in \cite{temnon} and \cite{tememb} in order to preserve the two above
properties as well. In particular, a desingularization $\calF$ (resp. embedded desingularization $\calE$) of all
generically reduced qe schemes over $\bfQ$ (resp. closed subschemes in regular qe schemes) is now available and
$\calF$ and $\calE$ go by blowing up regular centers and are functorial with respect to all regular morphisms.
The functoriality property is a serious achievement since it rigorously implies desingularization in many other
categories in characteristic zero, including Artin stacks, schemes acted on by regular group schemes, qe formal
schemes and complex or non-archimedean analytic spaces.

\subsection{Motivation and applications}
Non-functorial desingularization of qe schemes in \cite{temdes} only allowed to desingularize affine formal
schemes. In order to obtain a global desingularization result for formal schemes one had to construct a
desingularization which is functorial at least with respect to formal localizations. These are regular morphisms
of not finite type and it seems that the most natural way to ensure such functoriality is to achieve
functoriality with respect to all regular morphisms, as was done in \cite{temnon} and \cite{tememb}. Already
desingularization of formal varieties over $\bfC[[T]]$ is a new result, and currently it seems that
desingularization of formal varieties will be most useful for applications. Actually, it were few requests about
formal desingularization that convinced me to continue the research of \cite{temdes}. In particular, it seems
that desingularization of formal varieties may have applications to motivic integration (see \cite{Nic}), log
canonical thresholds (see \cite{FEM}), desingularization of meromorphic connections (see \cite{Ked}) and motivic
Donaldson-Thomas invariants studied by Kontsevich-Soibelman. Finally, it seems that the desingularization of
rigid spaces and Berkovich analytic spaces (not necessarily good) is also new.

\subsection{Overview}

We introduce all necessary terminology (e.g. qe schemes, blow ups, regular locus, etc.) in \S\ref{setup}. The
reader can look through this section and return to it when needed. In \S\ref{main} we formulate our main
results, explain how our method works in general and divide it to two stages. Then, both stages are studied in
details in \S\ref{algebra} and \S\ref{local}. In addition, we consider in \S\ref{examsec} few examples that
illustrate our algorithm. Finally, in \S\ref{apply} we deduce similar results for other categories including
stacks, formal scheme and various analytic spaces both in compact and non-compact settings.

\section{Setup}\label{setup}
Throughout this paper all schemes and formal schemes are assumed to be locally noetherian.

\subsection{Varieties}
{\em Variety} or {\em algebraic variety} in this paper always means a scheme $X$ which admits a finite type
morphism $X\to\Spec(k)$ to the spectrum of a field. If such a morphism is fixed then we say that $X$ is a {\em
$k$-variety} and $k$ is the {\em ground field} of $X$. It is an easy fact (see \cite{bmt}) that any reduced
connected variety $X$ possesses a maximal (and hence canonical) ground field $k\subset\calO_X(X)$.
Unfortunately, this is not true for non-reduced varieties as the following example shows.

\begin{exam}
Let $k$ be a field of characteristic zero with an irreducible curve $C$ which is not reduced at its generic
point $\eta$. Then $\Spec(\calO_{C,\eta})$ possesses various structures of a zero-dimensional $k(\eta)$-variety,
but none of them is "better" than another.
\end{exam}

\begin{rem}\label{groundrem}
(i) The above example extends to formal varieties. Moreover, even smooth formal varieties do not have to have a
canonical ground field. For example, already for $k=\bfQ(x)$ there exist many embeddings $k\into k[[t]]$, which
are as "good" as the obvious embedding, and, more generally, the field of coefficients in Cohen's theorem is not
unique.

(ii) As we will see later, the above observation is responsible for the main obstacle to proving functorial
desingularization by our method. In addition, it indicates that even for varieties it is more natural to study
absolute algorithms rather than the algorithms that take $k$ into account (for example by working with
$k$-derivations). We will return to this discussion in \S\ref{absolute}.
\end{rem}

\subsection{Regularity}

\subsubsection{Regular schemes}
There are many equivalent ways to say that a local ring $A$ is regular and here are two possibilities: (a) the
associated graded ring $\oplus_{n=0}^\infty m^n/m^{n+1}$ is isomorphic to $k[T_1\. T_d]$, where $m$ is the
maximal ideal and $k=A/m$ is the residue field, (b) the dimension of $A$ (i.e. the maximal length of a chain of
prime ideals decreased by one) equals to the dimension of the cotangent $k$-vector space $m/m^2$. We define the
{\em regular locus} $X_\reg$ of a scheme $X$ as the set of points $x\in X$ with regular $\calO_{X,x}$ and say
that $X$ is {\em regular} at each $x\in X_\reg$. The {\em singular locus} $X_\sing$ is defined as the complement
of $X_\reg$. Although regularity is an analog of smoothness, it is an absolute property, while smoothness is a
relative property. For example, a variety of positive characteristic can be smooth and not smooth over different
fields of definition.

\subsubsection{Monomial divisors}
A divisor $Z$ in a regular scheme $X$ is called {\em snc} (or strictly normal crossings) if its irreducible
components are regular and transversal, i.e. each $Z_{i_1}\cap Z_{i_2}\cap\dots\cap Z_{i_n}$ is regular of
codimension $n$ in $X$ (or empty). This is equivalent to the condition that Zariski locally at each point $x$
the divisor $Z$ is given by an equation $\prod_{i=1}^lt_i=0$ where $t_1\. t_n\in\calO_{X,x}$ is a regular family
of parameters. If $Z$ is given by an equation $\prod_{i=1}^nt_i^{m_i}=0$ then we say that it is {\em strictly
monomial}. Finally, if the above conditions hold only \'etale-locally (i.e. the parameters can be chosen in the
strict henselization $\calO^\sh_{X,x}$) then we say that $Z$ is {\em normal crossings} or {\em monomial},
respectively. Note that a closed subscheme $Z\into X$ is a (strictly) monomial divisor if and only if it is a
Cartier divisor whose reduction is (strictly) normal crossing.

\subsubsection{Regular morphisms}
A morphism $f:Y\to X$ is called {\em regular} if it is flat and has geometrically regular fibers. Since a finite
type morphism is regular iff it is smooth, this can be viewed as a generalization of smoothness to "large"
morphisms. A homomorphism of algebras $f:A\to B$ is {\em regular} if $\Spec(f)$ is regular (Hironaka calls
regular morphisms "universally regular", but our terminology is the standard one). It is well known that
regular/singular locus is compatible with regular morphisms, i.e. for a regular morphism $f:Y\to X$ we have that
$Y_\sing=f^{-1}(X_\sing)$ and $Y_\reg=f^{-1}(X_\reg)$. Similarly, one shows that the monomiality locus of a
divisor is compatible with regular morphisms. We warn the reader that the same is not true for strictly monomial
locus, since the preimage of a not strictly monomial divisor under an \'etale morphism can be strictly monomial.

\subsubsection{Equisingularity}\label{equi}
We say that a scheme $X$ is {\em equisingular} at a point $x$ if its reduction $X_0$ is regular at $x$ and
$X_\red$ is normally flat along $X_\red$ at $x$. Recall that the latter means that the $\calO_{X_\red}$-sheaves
$\calN_X^{i}/\calN_X^{i+1}$ are locally free at $x$, where $\calN_X\subset\calO_X$ is the radical. The set of
all points $x\in X$ at which $X$ is equisingular will be called the {\em equisingular locus} of $X$.
Equisingular loci behave similarly to regular loci. In particular, they are compatible with regular morphisms,
etc.

\begin{rem}
We prefer the notion "equisingular" since it is much shorter than other alternatives. Also, it reflects the
geometric meaning pretty well because in some sense the singularity of an equisingular scheme along any
irreducible component is constant; that is, the singularity on the entire component is as bad as at its generic
point. Note also that an equisingular scheme is regular if and only if it is generically reduced.
\end{rem}

\subsection{Quasi-excellent schemes}
\subsubsection{The definition}
For shortness, we will abbreviate the word quasi-excellent as {\em qe}. Quasi-excellent schemes were introduced
by Grothendieck in \cite[$\rm IV_2$, \S7.9]{ega} though the word "quasi-excellent" was invented later. These are
locally noetherian schemes $X$ satisfying the following conditions N and G (after Nagata and Grothendieck): (N)
for any $Y$ of finite type over $X$ the regular locus $Y_\reg$ is open, (G) for any point $x\in X$ the
completion homomorphism $\calO_{X,x}\to\hatcalO_{X,x}$ is regular. A qe scheme which is universally catenary
(see \S\ref{dimsec}) is called excellent.

\subsubsection{Connection to desingularization}
Obviously, the condition (N) is necessary in order to have a universal desingularization theory over $X$ (i.e.
in order to be able to desingularize schemes of finite type over $X$). Grothendieck proved the same for the
condition (G) in \cite[$\rm IV_2$, 7.9.5]{ega}: if any integral scheme of finite type over $X$ admits a regular
modification then $X$ is qe. It was suggested by Grothendieck and is believed by many mathematicians that the
converse is also true. Moreover, it is a common belief (or at least hope) that qe schemes admit much stronger
variants of desingularization discussed in \S\ref{desingsec}.

\subsubsection{Basic properties}\label{basic}
Main properties of quasi-excellence and excellence are as follows:

(1) They are invariant under many operations including passing to a scheme of finite type, localization and
henselization along a closed subscheme.

(2) If a ring $A$ is qe (i.e. $\Spec(A)$ is qe) and $\hatA$ is its completion along any ideal then the
completion homomorphism $A\to\hatA$ is regular.

(3) It is a very difficult result recently proved by Gabber that a noetherian $I$-adic ring $A$ is qe iff $A/I$
is qe. In particular, quasi-excellence is preserved under formal completions. See Remark \ref{gabrem} below for
more details.

(4) Basic examples of excellent rings are $\bfZ$, fields, noetherian convergent power series rings in complex
and non-archimedean analytic geometries, and schemes obtained from those by use of operations (1) and (2).

\begin{rem}\label{floor}
(i) Intuitively, general qe schemes have no "floor" unlike the algebraic varieties. For example, one cannot
fiber them by curves and many pathologies can occur with the dimension, as we will see below. Usually, in the
study of qe schemes one uses that they are in a "good relation with their roof" by the G-condition; that is,
their formal completions are regular over them in the affine case (i.e. the homomorphism $A\to\hatA$ is regular
when $X=\Spec(A)$). For example, the completion of a qe scheme $X$ along a subvariety (e.g. a closed point) is a
formal variety $\gtX$, and the desingularization theory for $X$ is closely related to that of $\gtX$ because of
the G-condition.

(ii) Formal varieties, in their turn, can be studied by various methods. In particular, one can fiber them by
formal curves (Gabber's adoption of de Jong's approach), one can algebraize them in the rig-regular case (our
adoption of Elkik's theory), and, very probably, one can generalize for them the algorithms for varieties by
switching to the sheaves of continuous derivations.
\end{rem}

\subsubsection{Caveats with the dimension theory }\label{dimsec}
A scheme $X$ is called {\em catenary} if for any point $x\in X$ with a specialization $y$ all maximal chains of
specializations between $x$ and $y$ have the same length. A scheme $X$ is {\em universally catenary} if any
scheme of finite type over $X$ is catenary. Actually, it is the catenarity condition which makes dimension
theory reasonable. The following simple example from \cite[$\rm IV_2$]{ega} shows that non-catenary schemes can
be not as horrible as one might expect.

\begin{exam}
Let $k$ be a field with an isomorphism $\phi:k\toisom k(t)$ (so, $k$ is of infinite absolute transcendence
degree, and one can take $k=F(t_1,t_2\dots)$, where $F$ is any field). Let $z$ be a closed $k$-point in
$\bfA^n_k=\Spec(k[x_1\. x_n])$ and $y$ be the generic point of an affine line not containing $z$. Let $\tilX_n$
be a localization of $\bfA^n_k$ with $n\ge 2$ on which both $y$ and $z$ are closed points and let $X_n$ be
obtained from $\tilX_n$ by gluing $y$ and $z$ via $\phi$ (i.e. we consider only functions $f\in\calO_{\tilX_n}$
with $f(y)=\phi(f(z))$). Note that $X_n$ is a "nice" qe scheme; for example, its normalization is a localized
variety $\tilX_n$. However, our operation obviously destroys the dimension theory on $X_n$, and indeed one can
easily show that $X_n$ is not catenary for $n\ge 3$ and is catenary but not universally catenary for $n=2$.
\end{exam}

\begin{rem}
(i) The above example is in a sense the most general one. Namely, it follows from \S\ref{basic}(2) that a local
qe ring $A$ is normal if its completion $\hatA$ is normal, and one can use this to show that any normal qe
scheme is catenary. Thus, the only source of non-catenarity on qe schemes is that sometimes branches of
different codimension on the same irreducible component can be glued on non-normal schemes. In particular,
non-catenarity is close in nature to local non-equidimensionality.

(ii) If $A$ is normal (or even regular) but not qe then it can happen that $\hatA$ is not normal or even is not
reduced. Also, there are normal but non-catenary not qe schemes.
\end{rem}

Another danger with qe schemes is that even an excellent ring can be infinite dimensional (by famous Nagata's
example). In particular, one cannot argue by induction on dimension and should use noetherian induction or
induction on codimension instead.

\begin{rem}
The reader may wonder why these pathological examples are worth any discussion. I agree that the schemes from
the above examples are curious but seem to be absolutely useless. However, the necessity to have them in mind
seems to be very useful from my point of view. It makes one to argue correctly and allows to quickly reject
approaches that could work for varieties but will not work for qe schemes (including the reasonable ones). For
example, the main induction in our desingularization method will be by codimension. Also, the non-catenary
example indicates that one must be extremely careful when dealing with non-equidimensional morphisms and schemes
(including varieties!). We will discuss this caveat in \S\ref{nonequi}.
\end{rem}

\subsubsection{Caveats with derivatives and bad DVR's}
In all examples from \S\ref{basic}(4) one establishes excellence by constructing a good theory of derivatives
(algebraic or continuous). The latter does not have to exist on regular not qe schemes, and this can be
interpreted as non-existence of global tangent space -- the spaces $m_x/m_x^2$ do not glue to a nice sheaf. The
source of the problem is that although the cotangent sheaf $\Omega^1_X$ is always quasi-coherent, it can be very
large (e.g. $\Omega^1_{\bfC/\bfQ}$ is a $\bfC$-vector space of continual dimension), and then its dual sheaf
$\Der_{X}=\calHom_{\calO_X}(\Omega^1_X,\calO_X)$ can be arbitrarily bad (e.g. not quasi-coherent, or even a
non-zero sheaf in a neighborhood of a point $x$ but with zero stalk at $x$). Moreover, the following example
shows that this can happen already for a qe trait, which is a very innocently looking scheme. (Recall that a
trait $X$ is the spectrum of a DVR, that is equivalent to $X$ being regular, local and of dimension one.)

\begin{exam}\label{badexam}
Let $k$ be a field and let $y=\sum_{i=0}^\infty a_ix^i\in k[[x]]$ be an element transcendental over $k[x]$. Then
the field $K=k(x,y)$ embeds into $k((x))$ and $\calO:=k[[x]]\cap K$ is a DVR with fraction field $K$ and
completion $k[[x]]$.

(i) If $\cha(k)=p$ and $y\in k[[x^p]]$ then $k((x))=\hatK$ is not separable over $K$ because it contains
$y^{1/p}$. In particular, the generic fiber of the completion homomorphism $\calO\to k[[x]]$ is not
geometrically reduced. This proves that the DVR $\calO$ is not a qe ring.

(ii) Though one can show that $\calO$ is excellent when $\cha(k)=0$, it still can have nasty differentials. For
example, let us assume in addition that the derivative $y'=\sum_{i=1}^\infty ia_ix^i$ is transcendental over
$K$. Consider the elements $y_i=x^{-n}(y-\sum_{i=0}^{n-1} a_ix^i)\in\calO$. An easy computation shows that
$\calO=k[x,y_0,y_1,y_2,\dots]_{(x)}$ and $\Omega^1_{\calO/k}$ is the $\calO$-submodule of $\Omega^1_{K/k}\toisom
Kdx\oplus Kdy$ generated by $dx,dy_0,dy_1,\dots$. Then it follows that actually
$\Omega^1_{\calO/k}=\Omega^1_{K/k}$. (Note also that in the case when $y'\in K$, the same computation shows that
$\Omega^1_{\calO/k}$ is obtained from the free $\calO$-module with generators $dx$ and $dy$ by adjoining the
elements $x^{-n}(dy-y'dx)$ for all natural $n$.) In particular, we obtain that
$\Der_k(\calO,\calO)=\Hom_\calO(K^2,\calO)=0$ while $\Der_k(K,K)=\Hom_K(K^2,K)\toisom K^2$ (and one shows
similarly that $\Der_k(\calO,\calO)\toisom\calO$ when $y'\in K$). This shows that the sheaf of $k$-derivations
on $\Spec(\calO)$ is not quasi-coherent and even has zero stalk at the closed point of $\Spec(\calO)$. Thus,
there is no good theory of derivations on $\calO$.
\end{exam}

\begin{rem}
(i) We saw that $\Der_{X}$ can behave wildly even for a regular qe scheme $X$. Since all current
desingularization algorithms over fields are based on derivatives, it is not clear if they can be extended to
all qe schemes without serious modifications. On the other hand one might hope that they can be
straightforwardly generalized to schemes that admit a closed immersion into a regular qe scheme "with good
theory of derivations". See, \cite[1.3.1(iii)]{temnon} for a precise conjecture.

(ii) The ring $\calO$ from Example \ref{badexam}(i) is a simplest example of a non-excellent ring. In addition,
$\calO$ is a very naively looking ring -- it is a DVR, and so it is regular, local and of dimension one. As an
additional demonstration of wildness of $\calO$ we note that its normalization in $k(x,y^{1/p})$ is a DVR which
is integral but not finite over $\calO$.
\end{rem}

\subsection{Blow ups}

\subsubsection{Basics}
Basic facts about blow ups can be found in \cite[\S2.1]{temdes} or in the literature cited there. Recall that
the {\em blow up} $f:\Bl_V(X)\to X$ along a closed subscheme $V$ is the universal morphism such that
$V\times_XX'$ is a Cartier divisor. In particular, $f=\Id_X$ iff $V$ is a Cartier divisor, $\Bl_X(X)=\emptyset$,
and $f$ is an isomorphism over $X\setminus V$. Also, $X\setminus V$ is dense in $\Bl_V(X)$ and so $f$ is
birational if $V$ is nowhere dense in $X$. The blow up always exists and it is the projective morphism given by
$\Bl_\calI(X)\toisom\Proj(\oplus\calI^n)$ where $\calI\subset\calO_X$ is the ideal of $V$ (we use our convention
that $X$ is locally noetherian and so $\calI$ is locally finitely generated). Conversely, any projective
modification is a blow up if $X$ possesses an ample sheaf, and in any case, blow ups form a very large cofinal
family among all modifications of a scheme (though the center of a typical blow up is highly non-reduced). A
blow up $\Bl_{\emptyset}(X)\toisom X$ is called {\em empty} or {\em trivial}.

\begin{rem}
Even empty blow ups play important role in functorial desingularization -- they are responsible for
synchronization.
\end{rem}

In the sequel we adopt the convention of \cite{temnon} that a blow up of $X$ consists of a morphism
$f:\Bl_V(X)\to X$ and a center $V$, i.e. the blow up "remembers" its center. This approach is finer than the
approach of \cite{temdes}, where a blow up was defined as a morphism isomorphic to a morphism of the form
$\Bl_V(X)\to X$ for some choice of $V$. As one may expect, we will see that the first approach is much better
suited for studying functorial desingularization.

\subsubsection{Operations with blow ups}\label{pushsec}
Blow ups are compatible with flat morphisms $f:X'\to X$ in the sense that
$\Bl_V(X)\times_XX'\toisom\Bl_{V'}(X')$ where $V'=V\times_XX'$.

If $f:\Bl_V(X)\to X$ is a blow up and $Z\into X$ is a closed subscheme, then the scheme-theoretic preimage
$Z\times_X\Bl_V(X)$ is called {\em full} or {\em total transform}, and we will denote it as $f^*(Z)$.

If $f:\Bl_V(X)\to X$ is a blow up, $Z\into X$ is a closed subscheme and $Z\setminus V$ denotes the open
subscheme of $Z$ obtained by removing $V$, then $Z\setminus V$ lifts to a subscheme in $\Bl_V(X)$ and its
schematical closure is called the {\em strict transform} of $Z$ under $f$ and will be denoted $f^!(Z)$. For
example, $f^!(Z)=\emptyset$ iff $|Z|\subset|V|$. Strict transforms are compatible with blow ups in the sense
that $f^!(Z)\toisom\Bl_{V|_Z}(Z)$.

If $i:U\into X$ is a locally closed immersion then any blow up $f:\Bl_W(U)\to U$ can be canonically {\em pushed
forward} to a blow up $i_*(f):\Bl_V(X)\to X$ where $V$ is the schematic closure of $W$ in $X$ ($W$ can be not
reduced, so we must take the {\em schematical closure}, i.e. the minimal subscheme $V\into X$ with $V|_{U}=W$).
Note that such extension is canonical because the blow up remembers its center. The restriction of $i_*(f)$ over
$U$ is $f$ itself, i.e. $\Bl_W(U)\toisom\Bl_V(X)\times_XU$ -- this is obvious for open immersions and this
follows from the properties of the strict transforms for a general locally closed immersion. Note that even if
$f$ is an isomorphism (i.e. $U$ is a Cartier divisor) $i_*(f)$ does not have to be an isomorphism.

\subsubsection{Blow up sequences}
Although a composition of blow ups is known to be isomorphic to a blow up, it is not clear how to choose a
center in a canonical way. If one ignores the centers one can study desingularizations by a single blow up, as
it is done in \cite{temdes}. However, for the sake of a more explicit description of a desingularization one
usually keeps all centers, i.e. considers whole {\em blow up sequences} $X'=X_n\to\dots\to X_1\to X_0=X$ with
the centers $V_i\into X_i$. Usually we will use the notation $X'\longto X$ for blow up sequences. All operations
with blow ups described above can be generalized to the blow up sequences straightforwardly (just iterate the
construction step by step). We say that a blow up sequence $f=f_{n-1}\circ\dots\circ f_1\circ f_0$ is {\em
$Z$-supported} for a closed subset $Z\subset X$ if all centers lie over $Z$.

\subsection{Desingularization}\label{desingsec}

\subsubsection{Weak desingularization}
A {\em weak desingularization} of an integral scheme $X$ is a modification $f:X'\to X$ with regular source. If
in addition $Z\times_ XX'$ is monomial for a closed subscheme $Z\into X$ then we say that $f$ is a weak
desingularization of the pair $(X,Z)$.

\begin{rem}
(i) Weak desingularization suffices to characterize qe schemes via \cite[$\rm IV_2$, \S7.9.5]{ega}.

(ii) Weak desingularization of varieties of characteristic zero can be proved by direct induction on dimension.
One fibers $X$ by curves and uses semi-stable modification theorem of de Jong and toroidal quotients, see
\cite{AdJ}.

(iii) The essential weakness of weak desingularization is that it does not control the modification locus (i.e.
the set of points of $X$ over which the modification is not an isomorphism). In particular, a desingularization
which modifies $X_\reg$ cannot be canonical.

(iv) The same result makes sense for any reduced scheme, but this generalization is not interesting since we can
simply use normalization to separate the irreducible components.
\end{rem}

\subsubsection{Desingularization}
By a (non-embedded) {\em desingularization} of a generically reduced scheme $X$ we mean a blow up sequence
$f:X'\longto X$ with regular source and such that $f$ is $X_\sing$-supported. If all centers are regular then we
say that $f$ is a {\em strong desingularization}.

\begin{rem}
(i) Currently, all proofs of (not weak) desingularization go through embedding varieties into smooth ambient
varieties and establishing an embedded desingularization, see \S\ref{embed}. In its simplest form such approach
leads to a non-strong desingularization, see \cite[\S8.2]{bmfun}.

(ii) An additional strengthening of the notion of desingularization is to require that each blow up center is
contained in the Hilbert-Samuel stratum of the largest order, where we naturally normalize the Hilbert-Samuel
function by codimension (so that it becomes constant on regular schemes) and use the natural partial order on
the set of all such functions. The methods of Hironaka, Bierstone-Milman and Villamayor provide such stronger
desingularization, but currently it is not achieved for qe schemes.
\end{rem}

\subsubsection{Desingularization of pairs}
By a {\em desingularization} of a pair $(X,Z)$ we mean a $(Z\cup X_\sing)$-supported blow up sequence $X'\longto
X$ with regular $X'$ and monomial $Z'=Z\times_XX'$. Classically one splits desingularization of $(X,Z)$ to usual
desingularization $X''\to X$ of $X$ and subsequent embedded desingularization of $Z\times_XX''\into X''$, but
such splitting is not necessary and sometimes seems to be not natural, see \S\ref{nonembbound}.

\subsubsection{Non-reduced schemes}
As defined above, the desingularization is rather meaningless for generically non-reduced schemes since it just
kills the generically non-reduced components. In particular, it can be easily obtained from desingularization of
generically reduced schemes, and hence does not involve anything new. A "right" desingularization of such
schemes is making them equisingular and it is usually achieved in the framework of strong desingularization. In
particular, it was established in the works of Hironaka, Bierstone-Milman and Villamayor.

\subsubsection{Embedded desingularization}\label{embed}
Let $X$ be a generically reduced variety of characteristic zero. Excluding special cases (e.g. low dimension),
all known constructions of a desingularization $X'=X_n\longto X_0=X$ go by embedding $X$ into a regular {\em
ambient} variety $M$ and successive blowing up $M$ so that the strict transform of $X$ becomes regular. Various
embedded desingularization algorithms have many similar features which we only outline here.

(1) {\em The boundaries.} One has to take the history of the process into account, see for example
\cite[\S3.6]{Kol}. This is done by considering on each $M_i$ a {\em boundary} $E_i$, which is the accumulated
exceptional divisor of the blow up sequence $M_i\longto M_0$. More concretely, $E_i$ is an ordered set of
divisors on $M_i$, which are called {\em components} of $E_i$ and are numbered by the history function. The
$i$-th boundary consists of the componentwise strict transform of the $(i-1)$-th boundary and the exceptional
divisor of the blow up $M_i\to M_{i-1}$. The pair $(M_i,E_i)$ is called the ambient variety with boundary, and
the basic objects of the embedded desingularization are the triples $(M_i,E_i,X_i)$. In classical embedded
desingularization the boundary is always {\em snc}, that is, its components are regular and meet transversally.
In applications one starts with $E_0=\emptyset$ but any choice of an {\em snc} $E_0$ is fine (see below).

(2) {\em Permissible centers.} The center $V_i$ of the blow up $M_{i+1}\to M_i$ is permissible in the sense that
it is regular and has {\em simple normal crossings} with $E_i$, i.e. for any component $D\in E_i$ we have that
locally at each point $x\in V_i\cap D$ either $V_i$ is transversal to $D$ or is contained in $D$. This ensures
that each $M_i$ is regular and each $E_i$ is {\em snc}.

(3) {\em Principalization.} Probably, the main paradigm of embedded desingularization is to replace the
desingularization problem with a very close problem of principalization of the ideal $\calI_X\subset\calO_M$
corresponding to $X$. Instead of the strict transform, one studies a {\em principal (controlled or weak)
transform} of $\calI=\calI_X$ under a blow up $M'\to M$ along a permissible center. This transform is obtained
from the full transform $\calI\calO_{M'}$ by dividing by an appropriate exceptional divisor. The ultimate aim of
the principalization is to find a permissible blow up sequence $f:M'\longto M$ such that the principal transform
of $\calI$ is $\calO_{M'}$ and hence $\calI'=\calI\calO_M$ is an exceptional divisor. In particular, $f$ induces
a desingularization of the pair $(M,X)$. Embedded desingularization, is obtained from principalization by
omitting the blow ups along components of the strict transform of $X$.

\begin{rem}\label{embrem}
(i) The main advantage of the principalization is that it replaces a geometric problem with an algebraic one. In
particular, it is much easier to compute principal transforms than the strict ones. In addition, all algorithms
deform $\calI$ severely in the process of principalization. This is done so that the ideal is replaced by an
equivalent one, which has the same principalization. No geometric interpretation of this procedure is known so
far. For example, $X$ with an isolated singularity $a$ is usually replaced with a highly non-reduced subscheme
of $M$ supported at $a$.

(ii) There are qe schemes that cannot be embedded into regular schemes (e.g. any non-catenary qe scheme). For
this reason one should separately establish non-embedded and embedded desingularization of qe schemes. The first
task was accomplished in \cite{temnon}. There are partial results on the embedded desingularization of qe
schemes in \cite{tememb}. In particular, the centers are regular but not transversal to the boundary.
Nevertheless, functorial desingularization of pairs is proved in \cite{tememb}.
\end{rem}

\subsubsection{Non-embedded desingularization with boundary}\label{nonembbound}
For the sake of completeness we note that one can generalize boundaries to the non-embedded setting. A {\em
boundary} $E$ on a scheme $X$ is an ordered set of locally principal closed subschemes of $X$ (with possible
repetitions). A finer form of desingularization of pairs is a desingularization of schemes with boundaries. The
latter seems to be a very recent notion, which was studied only in \cite{cjs} (for qe surfaces of all
characteristics) and \cite{tememb} (for qe schemes over $\bfQ$). To argue why non-embedded desingularization
with boundary might be a natural object to study we note that the embedded desingularization with boundary of
$(M,E,X)$ induces a non-embedded desingularization with boundary of $(X,E|_X)$ rather than just non-embedded
desingularization of $X$.

\subsubsection{Functorial desingularization}\label{functorial}
For the sake of concreteness we consider the non-embedded case in \S\S\ref{functorial}--\ref{absolute}. If
$f:Y\to X$ is a regular morphism then desingularizations $g':Y'\to Y$ and $g:X'\to X$ are {\em compatible} with
respect to $f$ if $g'$ is obtained from $g\times_XY$ by skipping empty blow ups and, moreover, $g'=g\times_XX'$
whenever $f$ is surjective (so, we even take the empty blow ups into account). If $\gtC$ is a class of schemes
(e.g. varieties, or qe schemes of characteristic zero) then by a {\em functorial desingularization on $\gtC$} we
mean a rule $\calF$ which to each $X\in\gtC$ assigns a desingularization $\calF(X):\oX\longto X$ in a way
compatible with all regular morphism between schemes from $\gtC$, i.e. for any such morphism $f:X'\to X$ the
desingularizations $\calF(X)$ and $\calF(X')$ are compatible with respect to $f$.

\begin{rem}
(i) Functoriality is a very strong property. It automatically implies desingularization in other categories
including equivariant desingularization, desingularization of stacks, analytic spaces, etc. Moreover, in a
seemingly paradoxical way it is usually easier to prove functorial desingularization since there is no problems
with gluing local desingularizations.

(ii) When proving functorial desingularization one must be very careful with synchronizing various blow ups. For
example, to construct $\calF(X\sqcup Y)$ from $\calF(X)$ and $\calF(Y)$ we must compare the singularities of $X$
and $Y$ and decide which one is blown up earlier (or simultaneously). In other words, we amplify $\calF(X)$ and
$\calF(Y)$ with synchronizing empty blow ups and then combine them into $\calF(X\sqcup Y)$. This illustrates the
role of the empty blow ups and explains why we worried for them in the definition of compatibility. See also
\cite[Lem. 2.3.1, Rem. 2.3.2]{temnon}. In addition, it is shown in \cite[Rem 2.3.4]{temnon} how the idea of
synchronization allows to represent any functorial desingularization as an algorithm governed by a
desingularization invariant.

(iii) Since Hironaka's foundational work many improvements and simplifications were made, and one of the main
achievements is that one obtains functorial desingularization. We try to outline (to some extent) the history of
the subject in \S\ref{hissec} below. Here we only note that in the recent papers \cite{Wl}, \cite{Kol} and
\cite{bmfun} one establishes functorial desingularization of varieties over a fixed field $k$ of characteristic
zero. Due to our convention from \S\ref{functorial}, this amounts to compatibility with all regular
$k$-morphisms between $k$-varieties, which are precisely all smooth $k$-morphisms. As for the class of all
varieties of characteristic zero and all regular morphisms between them, in addition to smooth morphisms these
works only checked compatibility with the ground field extensions, i.e. with the regular morphisms of the form
$X\otimes_kl\to X$ for a field extension $l/k$. It seems that full functoriality for varieties was established
only in \cite{bmt}.
\end{rem}

\subsubsection{On the history of desingularization of varieties of zero characteristic}\label{hissec}
It is very difficult to present a complete history of the field. So, I will only describe three stages and will
not even try to give all credits (including the very important contributions by Zariski and Giraud). The
original Hironaka's proof in \cite{Hir} was purely existential. The proofs of the second generation started with
the works \cite{Vil} and \cite{Vil2} of Villamayor and \cite{bm} and \cite{BM} of Bierstone-Milman. The main
focus in these (and many further) works is on constructing a canonical iterative desingularization algorithm
(with history) controlled by an appropriate invariant (or a set of invariants). Canonicity of the algorithm was
mainly understood as the fact that the constructed desingularization of a $k$-variety depended only on that
variety and was compatible with open immersions, which simplified the proofs a lot. Note that the new methods
heavily relied on some ideas (but not results) of Hironaka from \cite{Hir} an \cite{hir}, and, in addition,
Villamayor used Hironaka's results on idealistic presentation of Hilbert-Samuel function to obtain strong
desingularization. Starting with the recent work \cite{Wl} of Wlodarczyk (who builds a self-contained non-strong
desingularization algorithm but, again, heavily relies on the ideas of his predecessors), the main accent
shifted to functoriality of the desingularization and to recursive description of the algorithm, sometimes
making it unnecessary to introduce an invariant. In particular, it was shown in \cite{bmfun} how the algorithm
of \cite{BM} can be rewritten in a recursive form, and it was checked that this algorithm is functorial with
respect to all equidimensional smooth morphisms. Probably, all known canonical desingularization algorithms
become functorial with respect to all smooth (or even regular) morphisms after minor adjustments (see, for
example, \S\ref{nonequi}(2)), but this was not checked for most of the algorithms yet.

\subsubsection{Absolute desingularization of varieties}\label{absolute}
Intuitively it is clear that the functorial desingularization of varieties should be of absolute nature in the
sense that a ground field $k$ should not be taken into account. On the other hand, all known algorithms make
extensive use of the sheaves of $k$-derivations $\Der_{X/k}$, and in principle this may be an obstacle. More
concretely, the embeddings (or infinite localizations) like $$\bfA^{n-1}_{\bfQ(x_n)}=\Spec(\bfQ(x_n)[x_1\.
x_{n-1}])\into\Spec(\bfQ[x_1\. x_n])=\bfA^n_\bfQ$$ may be incompatible with the corresponding embedded
desingularizations because we work with the $\bfQ(x_n)$-derivatives in the source and with all
$\bfQ$-derivatives in the image. For example, for an ideal $I\subset A=\bfQ(x_n)[x_1\. x_{n-1}]$ the derivative
ideals $\Der_{A/\bfQ}(I)$ and $\Der_{A/\bfQ(x_n)}(I)$ used for these desingularizations are often different
because the derivation along the "constant direction" $x_n=0$ has a non-trivial effect. For this reason, it is
not clear if all known algorithms are of absolute nature and are compatible with all regular morphisms (though,
probably they are).

It was checked in \cite{bmt} that the algorithm of Bierstone-Milman is of absolute nature. However, the proof
used strong properties of the algorithm which are not known for some other algorithms. Moreover, it was shown in
loc.cit. that the algorithm admits an absolute description if one replaces the $k$-derivations sheaves
$\Der_{X/k}$ with the quasi-coherent absolute derivations sheaves $\Der_{X/\bfQ}$ (which can be very large). On
the other hand, the following interesting result from \cite{bmt} shows that any existing algorithm can be used
to produce an absolute algorithm just by using only its "$\bfQ$-component". In other words, an absolute
algorithm for varieties is the same as an algorithm for $\bfQ$-varieties.

\begin{theor}
Any functorial desingularization for $\bfQ$-varieties extends uniquely to a functorial desingularization of all
varieties of characteristic zero, their localizations and henselizations.
\end{theor}

This slightly surprising theorem is rather simple. The main idea is to use the approximation theory from
\cite[$\rm IV_2$, \S8]{ega} to approximate arbitrary varieties and regular morphisms between them with
$\bfQ$-varieties and smooth morphisms between them. In general, such approximation is possible for any
noetherian scheme over $\bfQ$, but, obviously, this is useless. So, the main observation about approximation of
varieties was that each variety is a projective limit of $\bfQ$-varieties with {\em smooth} and affine
transition morphisms. The latter smoothness condition is very special, and it reduces the problem to a standard
juggling with references to \cite[$\rm IV_2$, \S8]{ega}. Note that localizations and henselizations of varieties
are also such special projective limits, and so we can treat them in the same theorem. It is an interesting
question if there are other natural schemes which can be represented as such limits.

\subsubsection{Caveats with non-equidimensional schemes and morphisms}\label{nonequi}
In the desingularization theory one should be very careful when dealing with non-equidimensional varieties and
morphisms. We illustrate this by two examples.

(1) Usually, functorial embedded desingularization of $X$ in $M$ essentially depends on $X$ and its codimension
in $M$. For example, the resolution of $(M,X)$ will run faster than that of $(\bfA^1_M,X)$ when we run them
simultaneously (i.e. desingularize the disjoint union). Actually, for the algorithm from \cite[Ch. 3]{Kol} one
can show that if $j:M\into M'$ is a closed immersion with regular $M'$ then $\calE(M',X)$ and $\calE(M,X)$
induce the same desingularization of $X$ only when $j$ is of constant codimension. Probably, the same is true
for many other algorithms.

(2) It was recently noted by O. Gabber that the algorithm of Bierstone-Milman in \cite{bmfun} is only functorial
with respect to equidimensional smooth morphisms. A simple modification in the algorithm proposed in
\cite[\S6.3]{bmt} made the algorithm functorial with respect to all smooth (and regular) morphisms. Actually,
one just adjusts the synchronization slightly (in a sense, one replaces synchronization by dimension with
synchronization by codimension).

\subsubsection{Strict desingularization and a caveat with non-monomial locus}
By a {\em strict} desingularization of a pair $(X,Z)$ in \cite{temdes} one means a desingularization $X'\to X$
that modifies only the non-monomial locus of $(X,Z)$. This definition seemed natural to me but it turned out
that it does not make much sense. A detailed analysis can be found in \cite[\S{A.1.3}]{tememb}. Here we only
note that functorial strict desingularization does not exist even for varieties, and the assertion of \cite[Th.
2.2.11]{temdes} should be corrected as explained in \cite[Rem. A.1.1]{tememb}.

\subsection{Formal schemes}

\subsubsection{Quasi-excellent formal schemes}
By a {\em formal variety} we mean a noetherian formal scheme $\gtX$ whose special fiber $\gtX_s$ (defined by the
maximal ideal of definition) is a variety. An important stage of our method is desingularization of formal
varieties of characteristic zero, so we will explain briefly how the desingularization setup extends to formal
varieties. Formal varieties are {\em excellent} by results of Valabrega, see \cite{Val}. That is, for any affine
formal variety $\Spf(A)$ the ring $A$ is excellent. Since everything applies to general qe formal schemes, we
will work in such larger generality.

\begin{rem}\label{gabrem}
(i) To have a reasonable theory of qe formal schemes (other than formal varieties) one has to invoke Gabber's
theorem from \S\ref{basic}(iii). Otherwise, one does not even know that quasi-excellence is preserved by formal
localizations. Also, it is Gabber's theorem that implies that (quasi-) excellence is preserved by formal
completion.

(ii) The main intermediate progress towards Gabber's theorem was done in the paper \cite{NN} by
Nishimura-Nishimura, where the same result was proved conditionally assuming weak resolution of singularities
for local qe schemes. In particular, this settled the case of characteristic zero by using Hironaka's theorem
(which covers local qe schemes). Alternatively, one can use the results of \cite{temdes} as the
desingularization input.

(iii) Gabber strengthened the proof of \cite{NN} so that desingularization of local qe schemes is replaced with
a regular cover in the topology generated by alterations and flat quasi-finite covers. This argument is outlined
in Gabber's letter to Laszlo. The existence of such a regular cover for any qe scheme is a subtle and important
result by Gabber whose written version will (hopefully) be available soon. Actually, it is the only
desingularization result established for all qe schemes.
\end{rem}

\subsubsection{Regularity for qe schemes}
The underlying topological space of a formal scheme is too small to hold enough information even about reduced
formal subschemes. For this reason we define the singular locus of a formal scheme $\gtX$ as a closed formal
subscheme rather than as a subset in $|\gtX|$ (in particular, no regular locus is defined, though we remark for
the sake of completeness that one could work set-theoretically at cost of considering also a generic fiber of
$\gtX$ in one of the non-archimedean geometries). If $\gtX=\Spf(A)$ then we take for the singular (resp.
non-reduced or non-equisingular) locus the ideal defining $\Spec(A)_\sing$ (or other loci), and it turns out
that for qe formal schemes such definition is compatible with formal localizations and hence globalizes to
general qe formal schemes. Obviously, we use here that formal localization morphisms are regular on qe formal
schemes. Most probably, regularity and even reducedness does not make sense for general noetherian formal
schemes. We say that $\gtX$ is {\em regular} (resp. {\em reduced} or {\em equisingular}) if the singular (resp.
non-reduced or non-equisingular) locus is empty. We say that $\gtX$ is {\em rig-regular} if the singular locus
is given by an open ideal, and hence is a usual scheme. Intuitively, the latter means that the generic fiber of
$\gtX$ is regular (and this makes precise sense in non-archimedean geometry).

Regular and reduced loci are preserved by formal completions. Also, one uses a similar definition to introduce
the notion of regular morphisms between qe formal schemes and shows that the regular and reduced loci are
compatible with regular morphisms similarly to the case of schemes.

\subsubsection{Blow up sequences}
The notion of the formal blow up $\hatBl_\gtV(\gtX)$ along a closed formal subscheme $\gtV$ can be defined
similarly to the case of schemes. Then the formal blow up sequences are defined obviously. These notions are
compatible with formal completions, i.e. the $\calI$-adic completion of $\Bl_V(X)$ is $\hatBl_\gtV(\gtX)$, where
$\gtV$ and $\gtX$ are the $\calI$-adic completions of $V$ and $X$, respectively. All properties of usual blow
ups are generalized straightforwardly to the formal case, see \cite[\S2.1]{temdes}.

\subsubsection{Formal desingularization}
Since regular formal schemes and formal blow ups are defined, one defines desingularization of formal schemes
similarly to desingularization of schemes (including embedded desingularization, etc.).

\section{The method and the main results}\label{main}

\subsection{Results}

\subsubsection{The non-embedded case}

The main result of \cite{temnon} is that the class of all qe schemes over $\bfQ$ admits a strong non-embedded
desingularization. Here is a detailed formulation of this result.

\begin{theor}\label{mainth}
For any noetherian quasi-excellent generically reduced scheme $X=X_0$ over $\Spec(\bfQ)$ there exists a blow up
sequence $\calF(X):X_n\longto X_0$ such that the following conditions are satisfied:

(i) the centers of the blow ups are disjoint from the preimages of the regular locus $X_\reg$;

(ii) the centers of the blow ups are regular;

(iii) $X_n$ is regular;

(iv) the blow up sequence $\calF(X)$ is functorial with respect to all regular morphisms $X'\to X$, in the sense
that $\calF(X')$ is obtained from $\calF(X)\times_XX'$ by omitting all empty blow ups.
\end{theor}

\begin{rem}
An algorithm $\calF$ will be constructed from an algorithm $\calF_\Var$ for varieties, and we saw in
\S\ref{absolute} that $\calF_\Var$ is completely defined by its restriction $\calF_\bfQ$ onto the
$\bfQ$-varieties. So, in some sense $\calF$ is defined over $\bfQ$. Note, however, that $\calF$ is obtained by
"breaking $\calF_\Var$ to pieces" and reassembling them into a new algorithm, so it differs from $\calF_\bfQ$
even on $\bfQ$-varieties. This is necessary in order to have functoriality on all qe schemes.
\end{rem}

\subsubsection{The embedded case}
Here is the main result of \cite{tememb} formulated in the language of embedded desingularization, see \cite[Th.
1.1.6]{tememb}. Up to the non-embedded desingularization Theorem \ref{mainth}, this can be reformulated in the
language of non-embedded desingularization with boundary. We do not discuss such approach in this survey and
refer to \cite{tememb} for details.

\begin{theor}\label{embth}
For any quasi-excellent regular noetherian scheme $X$ of characteristic zero with an snc boundary $E$ and a
closed subscheme $Z\into X$ there exists a blow up sequence $f=\calE(X,E,Z):X'\longto X$ such that

(i) $X'$ is regular, the new boundary $E'$ is snc and the strict transform $Z'=f^!(Z)$ is regular and has simple
normal crossings with $E'$,

(ii) each center of $f$ is regular and for any point $x$ of its image in $X$ either $Z$ is not regular at $x$ or
$Z$ has not simple normal crossings with $E$ at $x$,

(iii) $\calE$ is functorial in exact regular morphisms; that is, given a regular morphism $g:Y\to X$ with
$D=E\times_XY$ and $T=Z\times_XY$, the blow up sequence $\calE(Y,D,T)$ is obtained from $g^*(\calE(X,E,Z))$ by
omitting all empty blow ups.
\end{theor}

\begin{rem}
The main weakness of this result is that the functor $\calE$ does not possess two important properties satisfied
by classical embedded desingularization functors.

(1) The centers of $\calE$ do not have to have normal crossings with the intermediate boundaries. In particular,
intermediate boundaries can be not snc, and even the iterative definition of these boundaries given in
\S\ref{embed}(2) should be corrected by replacing strict transform with principal transform (see
\cite[\S2.2]{tememb}).

(2) $\calE$ does not resolve the principal transform of $Z$. In particular, this cannot be used to obtain a
classical principalization of $\calI_Z$ as defined in \S\ref{embed}(3). However, if $Z$ is a Cartier divisor
(this situation is classically called "the hypersurface case") then $\calE$ induces a principalization.
\end{rem}

\subsubsection{Desingularization of pairs}
The strong principalization from \S\ref{embed}(3) is not achieved for qe schemes so far. However, the functors
$\calF$ and $\calE$ can be used to obtain a functorial desingularization of pairs.

\begin{theor} \label{princth}
For any quasi-excellent noetherian generically reduced scheme $X$ of characteristic zero with a closed subscheme
$Z\into X$ there exists a $(Z\cup X_\sing)$-supported blow up sequence $\calP(X,Z):X'\longto X$ such that $X'$
is regular, $Z\times_XX'$ is strictly monomial and $\calP$ is functorial in exact regular morphisms.
\end{theor}

The proof is very simple. First we blow up $X$ along $Z$ achieving that the full transform of $Z$ becomes a
Cartier divisor (this is an obvious principalization). Set $X'=\Bl_Z(X)$ and $Z'=Z\times_XX'$. Then we apply
$\calF$ to desingularize $X'$. Note that this step is needed even if we started with regular $X$. Let
$\calF(X'):X''\longto X'$ and $Z''=Z\times_XX''$. Finally, we apply $\calE(X'',\emptyset,Z'')$ to monomialize
$Z$. Note that a non-functorial desingularization of pairs is the main result of \cite{temdes}, and Theorem
\ref{princth} is a major strengthening of that result which was proved in \cite{tememb}.

\begin{rem}
The main disadvantage of our construction is that even when $X$ is regular $\calP$ can blow up a non-regular
center at the first step. This is only needed when $Z$ is not a Cartier divisor.
\end{rem}

\subsubsection{Semi-stable reduction}
In this section we just repeat the arguments from \cite[Ch. II, \S3]{KKMS}.  Assume that $\calO$ is an excellent
DVR of characteristic zero and $S=\Spec(\calO)$. Let $\eta$ and $s$ be its generic and closed points,
respectively. Assume that $X$ is a reduced flat $S$-scheme of finite type and with smooth generic fiber
$X_\eta$. It is well known that using a desingularization $Y\to X$ of the pair $(X,X_s)$ one can construct a DVR
$\calO'$ with a quasi-finite morphism $S'=\Spec(\calO')\to S$ and a modification $X'\to X\times_SS'$ such that
the morphism $X'\to S'$ is {\em semi-stable} (i.e. \'etale-locally $X'$ is of the form $\Spec(\calO'[t_1\.
t_n]/(t_1\dots t_m-\pi'))$ for a non-zero $\pi'\in\calO'$). Indeed, due to the assumption on the characteristic,
$Y$ is \'etale-locally of the form $\Spec(\calO[t_1\. t_n]/(t_1^{e_1}\dots t_m^{e_m}-\pi))$ for a non-zero
$\pi\in\calO$. Hence one can choose any $\calO'$ whose ramification degree $e$ over $\calO$ is divided by all
$e_i$ and take $X'$ to be the normalization of $Y\times_SS'$. Note that $X'$ as above is regular if and only if
$\pi'$ is a uniformizer. A complicated but purely combinatorial method to achieve that $X'$ is also regular is
described in \cite{KKMS}. It involves few blow ups along the strata of the reduction of $X'_s$ and its preimages
(note that they are snc divisors) and an additional extension of the DVR. The algorithm is described in terms of
the simplicial complex formed by the strata and the multiplicities of these strata in the closed fiber. In
particular, although originally formulated in the context of varieties, it applies to our situation verbatim.
This establishes the following theorem.

\begin{theor}
Assume that $\calO$ is an excellent DVR of characteristic zero, $S=\Spec(\calO)$ and $X$ is an $S$-scheme of
finite type and with smooth generic fiber. Then there exists a DVR $\calO'$ with a quasi-finite morphism
$S'=\Spec(\calO')\to S$ and a modification $X'\to X\times_SS'$ such that $X'$ is regular and the special fiber
$X_{s'}$ is an snc divisor (in particular, $X'$ is semi-stable over $S'$).
\end{theor}

\begin{rem}
The first step in our construction used $\calP$, so it is functorial. Functoriality of the whole construction
depends, thereby, only on the functoriality of the combinatorial algorithm. The algorithm from \cite[Ch.
III]{KKMS} seems to be not functorial (or canonical), but it seems very probable that functorial algorithms for
this problem should exist. So, I expect that the ramification degree of $\calO'/\calO$ and the modification
$X'\to X\times_SS'$ (for a fixed $S'$ with correct ramification) can be chosen functorially.
\end{rem}

\subsection{The method}
A very general idea of desingularizing qe schemes was discussed in Remark \ref{floor}: one wants to pass to
formal varieties by completion along subvarieties and desingularize the obtained formal varieties either by
algebraization or by generalizing the algorithms for algebraic varieties. The technical background is provided
by the following easy lemma.

\begin{lem}\label{lem1}
Let $X$ be a qe scheme such that $X_\sing$ is contained in a closed subvariety $Z\into X$ (e.g. $Z=X_\sing$) and
let $\gtX$ be the formal completion of $X$ along $Z$. Then $\gtX$ is a rig-regular formal variety and the formal
completion induces a bijective correspondence between desingularizations of $X$ and $\gtX$.
\end{lem}

The main point of the proof is that $\gtX_\sing$ is given by an open ideal and hence any desingularization of
$\gtX$ blows up only open ideals. Since open ideals live on a nilpotent neighborhood of $\gtX_s$ they algebraize
to closed subschemes of $X$ and hence the entire blow up sequence algebraizes as well. The lemma (and few more
similar claims) implies that any functorial desingularization $\hatcalF_\Var$ of formal varieties algebraizes
uniquely to a functorial desingularization $\calF_\small$ on the class of all qe schemes such that their
singular locus is a variety. Now we can explain very generally what are the two main stages of our method, and
we will describe them in more details in \S\S\ref{algebra}--\ref{local}.

Stage 1. {\it Algebraization.} The aim of this stage is to extend $\calF_\Var$ to $\calF_\small$. As we
explained above this reduces to constructing a desingularization functor $\hatcalF_\Var$ for rig-regular formal
varieties. The main tool is Elkik's theory which provides an algebraization of affine rig-regular formal schemes
with a principal ideal of definition. It allows to easily extend desingularization of varieties to rig-regular
formal varieties in a non-canonical way, see \cite[Th. 3.4.1]{temdes}. In \cite{temnon}, much more delicate
arguments were used in order to make this construction partially functorial. For technical reasons related with
Elkik's theory, $\hatcalF_\Var$ was only constructed for formal varieties with a fixed invertible ideal of
definition and for regular morphisms that respect these ideals.

\begin{rem}\label{functorrem}
(i) It is difficult to control functoriality since the algebraization procedure is absolutely non-canonical.
Much worse, as we observed in Remark \ref{groundrem}(i) even the ground field of algebraization is not
canonical. The latter turns out to be the main trouble since it is not easy to show that the existing
desingularization algorithms for varieties essentially depend only on the formal completion viewed as an
abstract formal scheme (i.e. without fixed morphism to a ground field). To illustrate the problem we note that
our formal algorithm must be equivariant also with respect to the automorphisms not preserving any ground field,
and such automorphisms do not have to be algebraizable by \'etale morphisms of varieties (compare with
\cite[3.56]{Kol} where all morphisms are defined over some ground field).

(ii) Most probably, $\hatcalF_\Var$ is fully functorial. The main problem of this stage is in establishing the
properties of $\hatcalF_\Var$ rather than in constructing it.
\end{rem}

Stage 2. {\it Localization} The aim of this stage is to reduce the general case to the case of schemes whose
singular locus is a variety, that is, to construct a functorial desingularization $\calF$ using a
desingularization functor $\calF_\small$ as an input. In a sense, we localize the desingularization problem at
this stage. Although, we cannot reduce to the case of a local scheme with an isolated singularity, we will only
use $\calF_\small(X)$ for rather special schemes $X$ such that $X_\sing$ is a variety. Namely, it will be enough
to know $\calF_\small(X)$ for a scheme $X$ which can be represented as a blow up of a local scheme $Y$ such that
$X_\sing$ is contained in the preimage of the closed point of $Y$ (and hence $X_\sing$ is a variety).

\begin{rem}\label{localrem}
The algorithm $\calF_\small$ is an extension of $\calF_\Var$, i.e. both agree on varieties. During the
localization stage a new algorithm $\calF$ is produced from $\calF_\small$. We will see in \S\ref{examsec} that
$\calF_\Var$ and $\calF$ differ already on algebraic curves. The construction of $\calF_\Var$ uses derivatives
and embedded desingularization, so it seems that it cannot be generalized to all qe schemes. On the other hand,
we have to build $\calF$ for all qe schemes in a "uniform way". Thus, it seems almost unavoidable that $\calF$
differs from $\calF_\Var$ on varieties.
\end{rem}

\section{Algebraization}\label{algebra}

Unless said to the contrary, we assume until the end of the paper that the characteristic is zero, i.e. all
schemes are $\bfQ$-schemes. The algebraization stage is rather subtle and technical and it is the bottleneck of
the method. In particular, it is "responsible" for most of the cases that elude from our method, including
generically non-reduced varieties, etc.

\subsection{Non-embedded rig-regular case}\label{nonemb}

Elkik's Theorem \cite[Th. 7]{Elk} implies that rig-regular formal varieties of characteristic zero with an
invertible (or locally principal) ideal of definition are {\em locally algebraizable} in the sense that they are
locally isomorphic to completions of varieties. For this reason we consider the pairs $(\gtX,\gtI)$ where $\gtX$
is a rig-regular formal variety and $\gtI$ is an invertible ideal of definition. We will be only able to
construct a formal desingularization $\hatcalF_\Var(\gtX,\gtI):\gtX_n\longto\gtX$ associated to such a pair and
functorial with respect to regular morphisms $\gtX'\to\gtX$ such that $\gtI'=\gtI\calO_{\gtX'}$. Most probably,
$\hatcalF_\Var$ is independent of $\gtI$ and is fully functorial, but this was not proved. For shortness, let us
say that $\gtX$ is a {\em principal formal variety} if it is rig-regular, affine, and with fixed principal ideal
of definition. Morphisms between such objects must be compatible with the fixed ideals.

Since we are going to establish functorial desingularization, it is enough to work locally. So, we can assume
that $\gtX$ is principal and then $\gtX$ is algebriazable by \cite[Th. 7]{Elk} and \cite[3.3.1]{temdes}, in the
sense that $\gtX=\hatX$ and $\gtI=\hatcalI$ for an affine variety $X$ with a principal ideal $\calI$. In
particular, the desingularization $\calF_\Var(X)$ induces a desingularization $\hatcalF_\Var(\gtX,\gtI)$ of
$\gtX$. The only thing we should do is to check that $\hatcalF_\Var(\gtX,\gtI)$ is well defined (i.e. is
independent of the choice of the algebraization) and functorial. The main idea beyond the argument is that all
information about $\gtX$ can be read off already from an infinitesimal neighborhood
$X_n:=(\gtX,\calO_\gtX/\gtI^n)=\Spec(\calO_X/\calI^n)$ with large $n$. This is based on \cite[3.2.1]{temnon}
which is an easy corollary of Elkik's theory. Roughly speaking, this result states that if $\gtX$ and $\gtX'$
are principal formal varieties and $n=n(\gtX)$ is sufficiently large then any isomorphism $X'_n\toisom X_n$
lifts to an isomorphism $\gtX'\toisom\gtX$.

Thus, it is clear that all information about the desingularization of $\gtX$ should be contained in some $X_n$,
though it is not so easy to technically describe this; especially because we want to prove functoriality of the
entire blow up sequence but only the first center is contained in $X_n$. We refer to \cite[\S3.2]{temnon} for a
realization of this plan. To illustrate some technical problems that one has to solve we note that if $X$ is an
algebraization of $\gtX$ and $f:X^{(p)}\longto X^{(0)}=X$ is its desingularization then the sequence of $n$-th
fibers $f_n:X^{(p)}_n\longto X^{(0)}_n$ is not determined by $X^{(0)}_n=X_n$ (for any $n$). However, one can
show that for sufficiently large numbers $k,n$ with $n\gg kp$ the tower
$$\calF_{n,k}(X^{(0)}_n,\calI_n):X^{(p)}_{n-kp}\to\dots\to X^{(1)}_{n-k}\to X^{(0)}_n$$ is uniquely determined
(up to a unique isomorphism) only by $X^{(0)}_n$ with the ideal $\calI_n=\calI\calO_{X^{(0)}_n}$ and, moreover,
is functorial in $X^{(0)}_n$ with respect to all regular morphisms. The above $\calF_{n,k}$ is a functor of
sequences of morphisms (not blow ups!) on certain non-reduced schemes with fixed principal ideal which are
called Elkik fibers in \cite{temnon}, and $\calF_{n,k}$ is the heart of the technical proof that $\hatcalF_\Var$
is a well defined functor.

\subsection{Limitations}
The limitations of our algebraization method are related to the assumptions in Elkik's theory. For example, in
the algebraization theorem \cite[Th. 7]{Elk} one assumes that the formal scheme is rig-smooth over the base and
possesses a principal ideal of definition. In addition, no result for algebraization of a pair $(X,Z)$, where
$Z\into X$ is a closed subscheme, is known. Let us discuss what is the impact of these assumptions on our
method.

\subsubsection{Closed subschemes} I do not know if algebraization of pairs is possible under reasonable
assumptions (say, $\gtX$ is regular and $\gtZ\into\gtX$ is rig-regular). I hope that some progress in this
direction is possible and this should be studied in the future. Currently, the lack of algebraization of pairs
is the main reason that our embedded desingularization theorem is much weaker than its classical analog. In
\cite{tememb} one only uses algebraization of pairs $(\gtX,\gtZ)$, where $\gtX$ is rig-regular and $\gtZ$ is
supported on the closed fiber (and hence is a scheme).

\subsubsection{Rig-regularity}
Simple examples show that rig-smoothness is necessary in order for algebraization to exist. Localization stage
reduces desingularization of generically reduced schemes to desingularization of rig-regular varieties, and
rig-smoothness is equivalent to rig-regularity in characteristic zero (see, for example, \cite[\S3.3]{temdes}).
Thus, the algebraization stage as it is works only in characteristic zero (even if resolution of varieties would
be known). Over $\bfQ$ the only limitation imposed by the rig-regularity assumption is that our method does not
treat generically non-reduced schemes. To deal with the latter case on should desingularize {\em
rig-equisingular} formal varieties (i.e. formal varieties whose non-equisingular locus is given by an open
ideal), but the latter formal schemes are not locally algebraizable in general.

\subsubsection{Principal ideal of definition}
I do not know if this assumption in Elkik's theorem is necessary. It causes to certain technical difficulties in
our proofs and extra-assumptions in intermediate results, but does not affect the final results. Step 1 in the
localization stage (see \S\ref{local}) is needed only because of this assumption. Also, it is this assumption
that makes us to define the functor $\hatcalF_\Var(\gtX,\gtI)$ rather than a functor $\hatcalF_\Var(\gtX)$.

\section{Localization}\label{local}
The localization stage is very robust, and can be adopted to work with almost all types of desingularization,
including embedded desingularization, and desingularization of non-reduced schemes. Also, it is not sensitive to
the characteristic. For simplicity, we will stick with the non-embedded case which is established in
\cite[\S4.3]{temnon}.

\subsection{Construction of $\calF$}\label{Fsec}
Consider the category $\gtC_\small$ whose elements are pairs $(X,D)$ where $X$ is a noetherian generically
reduced qe scheme, $D\into X$ is a closed subvariety which is a Cartier divisor, and morphisms $(X',D')\to(X,D)$
are regular morphisms $f:X'\to X$ such that $D'=f^*(D)$. By the algebraization stage and Lemma \ref{lem1}, the
original desingularization functor $\calF_\Var$ extends to a desingularization functor $\calF_\small$ on
$\gtC_\small$. The aim of the localization stage is to construct a desingularization $\calF$ of qe schemes using
$\calF_\small$ as an input.

The construction of $\calF$ goes by induction on codimension, i.e. we will construct inductively a sequence of
blow up sequence functors $\calF^d$ which desingularize $X$ over $X^{\le d}$, where the latter denotes the set
of points of $X$ of codimension at most $d$. Intuitively (and similarly to \S\ref{pushsec}), each $\calF^d(X)$
is the pushout of the desingularization $\calF(X)|_{X^{\le d}}$ under the embedding $X^{\le d}\into X$, i.e. it
is the portion of $\calF^d(X)$ defined by the situation over $X^{\le d}$. More specifically, each center of
$\calF^d(X)$ has a dense subset lying over $X^{\le d}$ and $\calF(X)$ is obtained from $\calF^d(X)$ by inserting
(in all places) few new blow ups whose centers lie over $X^{>d}:=X\setminus X^{\le d}$. Thus, the resulting
algorithm works as follows. Take empty $\calF^0(X)$. That is, start with $X$, which is the canonical
desingularization of itself over $X^{\le 0}$ (we use that $X$ is generically reduced by our assumption). First
we resolve the situation over the points of $X$ of codimension one by a functor $\calF^1$ (without caring for
other points). Then we improve $\calF^1$ over finitely many points of codimension two and leave the situation
over the codimension one points unchanged. This gives a functor $\calF^2$ which agrees with $\calF^1$ over the
codimension 1 points and resolves each generically reduced qe scheme $X$ over $X^{\le 2}$. We proceed similarly
ad infinitum, but for each noetherian $X$ the process stops after finitely many steps by noetherian induction.

Now let us describe how $\calF^d$ is constructed from $\calF^{d-1}$ and $\calF_\small$. Given a blow up sequence
$f:X'\longto X$ by its {\em unresolved locus} $f_\sing$ we mean the set of points of $X$ over which $f$ is not a
strong desingularization. In other words, $f_\sing$ is the union of the images of the singular loci of both $X'$
and the centers of $f$. By the induction assumption, $\calF^{d-1}(X)_\sing$ is of codimension at least $d$ and
hence it contains only finitely many points $x_1\. x_m$ of exact codimension $d$. We should only improve $\calF$
over these points, and the latter is done as follows.

Step 1. As a first blow up we insert the simultaneous blow up at all new points $x_1\. x_m$ (we act
simultaneously in order to ensure functoriality). Then the preimage of each $x_i$ on any intermediate blow up of
the sequence is a Cartier divisor (which will be needed later in order to use $\calF_\small$).

Step 2. Next, we improve all centers of $\calF^{d-1}(X)_\sing$ over $x_i$'s by resolving the singularities of
these centers over $x_i$'s. We use here that these singularities are of codimension at most $d-1$ in the centers
and so we can apply the functor $\calF^{d-1}$. To summarize, before blowing up each center $V_j$ of
$\calF^{d-1}$ we insert a blow up sequence which desingularizes $V_j$ over $x_i$'s.

Step 3. At the last step we obtain a sequence $X'\longto X$ of blow ups whose centers are regular over $X^{\le
d}$, but $X'$ may have singularities over $x_i$. Observe that the singular locus of the scheme
$X_{x_i}=\Spec(\calO_{X,x_i})\times_XX'$ is contained in the preimage of $x_i$, and the preimage of $x_i$ is a
Cartier divisor $E_i$ (thanks to Step 1) which is a variety over $k(x_i)$. So, $X_{x_i}$ can be resolved by a
blow up sequence $f_i=\calF_\small(X_{x_i},E_i)$. It remains to extend all $f_i$'s to a blow up of $X'$ and to
synchronically merge them into a single blow up sequence $X''\longto X'$. The composition $\calF^d(X):X''\longto
X$ is a required desingularization of $X$ over $X^{\le d}$ which coincides with $\calF^{d-1}$ over $X^{\le
d-1}$.

\begin{rem}
The center of a blow up is often reducible, and in Step 3 of the construction of $\calF^d$ we often obtain a
center with many components that are regular over $x_i$'s but probably have non-empty intersections over
$X^{>d}$. Thus, it is important that in Step 2 of the construction of $\calF^d$ we are able to desingularize the
reducible blow up centers inherited from $\calF^i$ with $i<d$. In particular, even if we are only interested to
desingularize integral schemes, we essentially use in our induction that the desingularization is constructed
for all reduced schemes.
\end{rem}

\subsection{Examples}\label{examsec}
We will compare $\calF$ and $\calF_\Var$ in the case of few simple varieties. For $\calF_\Var$ we take the
desingularization functor of Bierstone-Milman, which is functorial in all regular morphisms by \cite{bmt}.

\subsubsection{Plain curves}\label{curvsec}
Assume that $X$ is a generically reduced plain algebraic curve. A strong desingularization is uniquely defined
up to synchronization because one has to blow up the singular points until the curve becomes smooth. On the
other hand, synchronization of these blow ups depends on various choices. Note that the Hilbert-Samuel strata of
$X$ are the equimultiplicity strata because $X$ embeds into a smooth surface. It follows easily that
$\calF_\Var$ is synchronized by the multiplicity. Namely, one blows up the points of maximal multiplicity at
each step until all points are smooth. The synchronization of $\calF$ is slightly different. Because of Step 1
in the localization stage we simultaneously blow up all singular points once. After that we skip Step 2 and use
$\calF_\Var$ at Step 3. To summarize, we blow up all singularities once, and then switch to the synchronization
by multiplicity, similarly to $\calF_\Var$.

\begin{rem}
The same description holds true for varieties with isolated singularities. At the first stage, $\calF$ blows up
all these singularities. One obtains a blow up $X'\to X$ and then simply applies $\calF_\Var(X')$. Thus, $\calF$
blows up the same centers as $\calF_\Var$ but the synchronization can be different when $X$ has more than one
singular point.
\end{rem}

\subsubsection{Surfaces}\label{surfsec}
Let us consider examples of the next level of complexity. Namely, let $X$ be a surface such that $C=X_\sing$ is
a curve with the set of generic points $\eta=\{\eta_1\.\eta_n\}$. The functor $\calF^1$ acts as follows. On the
semi-local curve $X_\eta$, which is the semi-localization of $X$ at $\eta$, $\calF^1$ acts as was explained in
\S\ref{curvsec}. We extend the blow up sequence $\calF^1(X_\eta):X'_\eta\longto X_\eta$ to a blow up sequence
$\calF^1(X):X'\longto X$ in the natural way (that is, the centers of $\calF^1(X)$ are the Zariski closures of
the centers of $\calF^1(X_\eta)$). After that we produce $\calF(X)=\calF^2(X)$ by inserting new blow ups into
$\calF^1(X)$. This is done in three steps described in \S\ref{Fsec}. All new blow ups will be inserted over the
set $b=(b_1\. b_m)$ such that $\calF^1(X)$ is not a strong desingularization precisely over the points of $b$.
In particular, the first blow up is along $b$.

\begin{exam}
If $a\in X$ is a "generic point" of $C$ then $\calF(X)=\calF^1(X)$ over a neighborhood of $a$. For example, this
is the case of any surface of the form $X=\bfA^1_k\times Y$ for a curve $Y$. One easily sees that in this case
$\calF$ and $\calF_\Var$ differ only by synchronization, as in the case of curves.
\end{exam}

Next we consider two examples of a Cartier divisor in $\bfA^3_k=\Spec(k[x,y,z])$ with a non-isolated and
"non-generic" singularity at the origin $a$. For the sake of comparison, we will show how $\calF_\Var$ resolves
the same examples. For reader's convenience, a brief explanation of how $\calF_\Var$ can be computed in these
examples will be given in \S\ref{fsec}.

\begin{exam}\label{exam1}
Whitney umbrella $X$ is given by $y^2+xz^2=0$. In this case $C$ is the $x$-axis, and blowing it up resolves all
singularities. So, $\calF(X)=\calF^2(X)=\calF^1(X)$ just blows up $C$. However, $a$ is not a "generic point" of
the singular locus and other algorithms feel this. In particular, $\calF_\Var$ first blows up $a$. The blow up
$X_1=\Bl_a(X)$ is covered by two charts: the $x$-chart $X_{1x}$ and the $z$-chart $X_{1z}$. Since $X_{1x}$ is
defined by $y_1^2+x_1z_1^2=0$ for $x_1=x$, $y_1=y/x$ and $z_1=z/x$, we see that it has the same singularity as
$X$. Namely, the singular locus $C_1$ is the line $y_1=z_1=0$ (which is the strict transform of $C$). Since
$X_{1z}$ is defined by $y_1^2+x_1z_1=0$ for $x_1=x/z$, $y_1=y/z$, $z_1=z$ (for simplicity we denote the local
coordinates by the same letters as earlier), the singularity of $X_{1z}$ is the isolated orbifold point $c_1$
given by $x_1=y_1=z_1=0$. Because of a synchronization issue, $C_1$ is dealt with first (a non-monomial ideal
$\calN$ has order two along $C_1$ and order one at $c_1$, see \S\ref{fsec} for a similar computation). The
second blow up is, again, at the pinch point $a_1\in C_1$. The same computation as above shows that the singular
locus of $X_2$ consists of $C_2$, which is the strict transform of $C$ (and $C_1$), an orbifold point $c'_2$,
which sits over $a_1$ and an orbifold point $c_2$ which is the preimage of $c_1$. The third blow up is along
$C_2$, so the singular locus of $X_3$ consists of two orbifold points $c_3$ and $c'_3$ sitting over  $c_2$ and
$c'_2$, and the last blow up is along $\{c_3,c'_3\}$. This resolves $X$ completely. Note also that for $0\le
i\le 2$ the local structure of $X_i$ along $C_i$ is similar. Nevertheless, the algorithm blows up the pinch
point twice and then decides to blow up the whole singular line $C_2$ because of the history of the process.
\end{exam}

\begin{exam}\label{exam2}
Let $Y$ be given by $y^2+xz^3=0$. The entire resolution process is too messy, so we will only work with a
sequence of affine charts. As we saw in the previous example, there might be (and there are) blow ups on the
other charts that may (but do not have to) be simultaneous with the blow ups on our charts. In this case, we
also have that $C=Y_\sing$ is the $x$-axis and $\calF^1(Y):\tilY\to Y$ is the blow up along $C$. Let us describe
an {\em affine chart} $Y_n\longto Y_0=Y$ of $\calF(Y):\oY_\on\longto\oY_0=Y$. It is obtained by choosing each
time an appropriate affine chart of the blow up and restricting the remaining sequence over that chart. In
particular, $n<\on$ because of the synchronization with other charts, which we ignore for simplicity.

Note that $\tilY$ has an isolated orbifold singularity above $a$ locally given by $y_1^2+x_1z_1=0$, where
$x_1=x$, $z_1=z$ and $y_1=y/z$. Thus, $\calF(Y)\neq\calF^1(Y)$, $b=\{a\}$ and Step 1 inserts the blow up at $a$
as the first blow up. So, $\oY_1=\Bl_a(Y)$ and we will study how $\calF$ proceeds on the (most interesting)
$x$-chart $Y_1$ defined by the equation $y_1^2+x_1^2z_1^3=0$, where $x_1=x$, $y_1=y/x$ and $z_1=z/x$. Note that
the strict transform of $C$ is the $x_1$-axis, which we denote by $C_1$. Since $C_1$ is regular, no blow up is
inserted at Step 2, and so the second blow up is along $C_1$. Thus, $\Bl_{C_1}(Y_1)$ has only one chart $Y_2$
given by $y_2^2+x_2^2z_2=0$, where $x_2=x_1$, $y_2=y_1/z_1$ and $z_2=z_1$. In particular, $Y_2$ is a Whitney
umbrella, and $(Y_2)_\sing$ is contracted to the point $a$ by the projection $Y_2\to Y$. The singularity of
$Y_2$ is resolved at Step 3 by applying $\calF_\Var(Y_2)$, which is the same as $\calF_\Var(X)$ from Example
\ref{exam1}. So, the pinch point of $Y_2$ is blown up twice and the fifth blow up blows up the entire line,
which is the strict transform of $(Y_2)_\sing$ (note that this singular line appeared for the first time as the
preimage of $a$ under $Y_1\to Y$).

Now, let us describe $\calF_\Var(Y):Z_n\longto Z_0=Y$. The first two blow ups are at the origin (similarly, to
Whitney umbrella). Thus $Z_1=Y_1$ and for $Z_2$ we take the (most interesting) $x$-chart, so $Z_2$ is given by
$y_2^2+x_2^3z_2^3=0$, where $x_2=x_1$, $y_2=y_1/x_1$ and $z_2=z_1/x_1$. Note that the singularities of $Y_2$ and
$Z_2$ are different. The third blow up is along the line $y_2=z_2=0$ (which is the strict transform of the
original singular line $C$), hence its only chart is $Z_3$ given by $y_3^2+x_2^3z_2=0$, where $x_3=x_2$,
$y_3=y_2/z_2$ and $z_3=z_2$. The next blow up is along the line $x_3=y_3=0$ and $Z_4$ is given by
$y_4^2+x_4z_4=0$. Finally, the last blow up is at the point $x_4=y_4=z_4$.

We summarize by saying that (up to our matching of affine charts, which is sort of informal) $\calF$ blows up
point, old line (the strict transform of $C$), point, point, new line, while $\calF_\Var$ blows up point, point,
old line, new line, point.
\end{exam}

\begin{rem}
(i) We saw that when the singularities are not isolated, the two algorithms can blow up different centers.
Unlike the synchronization issues, this makes the algorithms very different because they deal with different
singularities after the first choice of different centers. In particular, it is unclear to me how one can
compare the algorithms in general (or match affine charts).

(ii) Although these two examples do not give enough intuition for deciding which algorithm is faster, I would
expect that in general $\calF$ produces more complicated desingularizations than $\calF_\Var$. Nevertheless, the
tendency of $\calF$ to blow up curves (and other higher dimensional centers) at first occasion can shorten the
desingularization process in some cases.
\end{rem}

\subsubsection{Computing $\calF_\Var$}\label{fsec}
Let $X$ and $Y$ be as in Examples \ref{exam1} and \ref{exam2}. We will compute $\calF_\Var(Y)$ and
$\calF_\Var(X)$ can be found similarly. We apply the algorithm from \cite[\S5]{bmfun}, so the reader is advised
to consult \cite{bmfun} for more details.  Since $Y$ is a hypersurface in $M=\bfA^3$, we can take the
equimultiple stratum of $Y\into M$ as a presentation of the maximal Hilbert-Samuel stratum of $Y$. Thus, it
suffices to resolve the marked ideal $(y^2+xz^3)$ with $d=2$ on $M$. We can take the $(xz)$ plane $P$ given by
$y=0$ as the hypersurface of maximal contact. Restricting the coefficient ideal
$\calC_\emptyset^{d-1}=(y^2+xz^3)+(2y,3zx^2,z^3)^2$ on $P$ gives the ideal $\calI=(z^6,xz^3)$ with $d_\calI=d=2$
(we refer to \cite{bmfun} and especially its \S3.4, \S4 and Case A of Step II in \S5). To resolve $(y^2+xz^3)$
on $M$ is the same as to resolve $\calI$ on $P$, and for the latter we go to Case B of Step II in
\cite[\S5]{bmfun}. The monomial part $\calM$ of $\calI$ is trivial (no history) and the non-monomial part is
$\calN=(z^6,xz^3)$ with $d_\calN=4$. So, the companion ideal $\calG$ is $\calN$ with $d_\calG=4$ and it easy to
see (via further maximal contact reduction) that $\calG$ is resolved by blowing up the origin. Hence the first
blow up in the resolution of $\calI$ is along the origin $a$ (and similarly for the resolution of $Y$). Consider
the $x$-chart $P_1$ of $\Bl_a(P)$. The principal transform of $\calI$ on this chart is
$\calI_1=(x_1^2z_1^3,x_1^4z_1^6)=(x_1^2z_1^3)$ (it is obtained from the full transform by dividing by
$z_1^d=z_1^2)$. Now, $\calI_1=\calM_1\calN_1$ for monomial $\calM_1=(x_1^2)$ and non-monomial $\calN_1=(z_1^3)$.
It follows that $\calG_1=(z_1^3)$ with the exceptional divisor $V(x_1)$. Hence $\calG_1$ is resolved by two blow
ups: first we blow up the point $x_1=z_1=0$ due to Case B of Step I in \cite[\S5]{bmfun} in order to separate
the singular locus $V(z_1)$ from the old boundary. This gives $\calG_2=(z_2^3)$ on the $x$-chart, and then we
blow up the line $V(z_2)$. Tracking the effect on the principal transforms of $\calI_1$ under these two blow ups
we see that $\calI_2=x_2^{-2}(x_2^5z_2^3)=(x_2^3z_2^3)$ and $\calI_3=(x_3^3z_3)$. The latter ideal is monomial
because the exceptional divisor is $V(x_3z_3)$ at this stage. So, Step II in \cite[\S5]{bmfun} deals with it by
Case A. It is easy to see that $\calI_3$ is resolved by two additional blow ups. First one blows up the line
$V(x_3)$, obtaining $\calI_4=(x_4z_4)$, and then one blows up the remaining singular point $V(x_4,z_4)$,
resolving $\calI$ completely. Thus, $\calI$ is locally resolved by an affine chart of a blow up sequence of
length five, and its centers are as follows: the point $V(x,z)$, the point $V(x_1,z_1)$, the line $V(z_2)$, the
line $V(x_3)$ and the point $V(x_4,z_4)$. Therefore the centers of the corresponding "affine chart" of
$\calF_\Var(Y)$ are as follows: the point $V(x,y,z)$, the point $V(x_1,y_1,z_1)$, the line $V(y_2,z_2)$, the
line $V(x_3,y_3)$ and the point $V(x_4,y_4,z_4)$.

\section{Desingularization in other categories}\label{apply}
We will show that the main Theorems \ref{mainth}, \ref{embth} and \ref{princth} imply analogous
desingularization theorems in many other categories, including qe stacks, formal schemes and various analytic
spaces in characteristic zero. Also, we will show that desingularization of non-compact objects follows as well.

\subsection{Stacks}
Let $\gtX$ be an Artin stack with a smooth atlas $p_{1,2}:R\rra U$, so $U\to\gtX$ is a smooth covering and
$R=U\times_\gtX U$. If $\gtV\into\gtX$ is a closed substack then we define the blow up $\Bl_\gtV(\gtX)$ using
the chart $\Bl_W(R)\rra\Bl_V(U)$ where $V=\gtV\times_\gtX U$ and $W=\gtV\times_\gtX R$ (we use here that the
blow ups are compatible with flat morphisms and so $\Bl_W(R)=\Bl_V(U)\times_{U,p_i}R$ for $i=1,2$). We say that
a stack is regular or qe if it admits a smooth cover by such a scheme. A strong desingularization is now defined
as in the case of schemes.

\begin{theor}
The blow up sequence functors $\calF$, $\calE$ and $\calP$ extend uniquely to noetherian qe stacks over $\bfQ$.
\end{theor}

To prove this theorem for $\calF$ we take any stack $\gtX$ as above and find its smooth atlas $R\rra U$. Then
$\calF(R)\rra\calF(U)$ is an atlas of a blow up sequence $\calF(\gtX)$. An interesting case is when $\gtX=X/G$
for an $S$-scheme $X$ acted on by an $S$-smooth group scheme $G$. Then the above theorem actually states that
$X$ admits a $G$-equivariant desingularization. Other functors are dealt with similarly.

\subsection{Formal schemes and analytic spaces}\label{formal}
In the categories of qe formal schemes, complex analytic spaces, rigid analytic spaces and analytic $k$-spaces
of Berkovich the notions of regularity and blow ups are defined. So, one can define strong desingularization
similarly to the case of schemes. Hironaka proved desingularization of complex analytic spaces, but this
required to insert major changes in his method (and the main reason is that his method is not canonical). The
new algorithms are known to work almost verbatim for complex analytic spaces, though strictly speaking, one
should repeat the entire proof word by word. The desingularization of affine formal schemes was deduced in
\cite{temdes} from non-functorial desingularization of affine schemes, but this approach did not yield global
desingularization of formal schemes.

It turns out that functorial desingularization of qe schemes is so strong that it rigorously implies functorial
desingularization of all above objects. The strategy is always the same, so let us stick with the non-embedded
desingularization. We cover a generically reduced object $X$ (i.e. the non-reduced locus is nowhere dense) by
compact local subobjects $X_1\. X_n$ (e.g. affinoid subdomains, affine formal subschemes or Stein compacts) and
observe that $A_i=\calO_X(X_i)$ are qe rings and for any smaller object $X_{ijk}\subset X_i\cap X_j$ the
localization homomorphisms $A_i\to \calO_X(X_{ijk})$ are regular (e.g. formal localization is regular on qe
formal schemes). Thus, completion/analytification of the desingularization $\calF(\Spec(A_i))$ yields a
desingularization of $X_i$, and these local desingularizations glue together because $\calF$ is compatible with
all regular morphisms.

\begin{rem}
(i) Recall that $\calF$ is of absolute nature, and actually it is constructed from a functor $\calF_\bfQ$ on
$\bfQ$-varieties. Thus, the obtained desingularization of all above objects is algebraic, and even defined over
$\bfQ$ in some sense. The latter might look surprising since there are non-algebraizable analytic singularities,
so we illustrate below the differences between our method and "naive algebraization".

(ii) First of all, thanks to the localization stage we only have to algebraize rather special classes of
singularities, which generalize in some sense the isolated singularities. Furthermore, even when $x\in X$ is an
isolated complex singularity, we do not algebraize a complex neighborhood of $x$ but only its formal
neighborhood $\hatX_x=\Spf(\hatcalO_{X,x})$. This operation is "too local" at $x$, so it does not have to extend
to an analytic neighborhood of $x$. Moreover, for the sake of functoriality we had to study all algebraizations
of $\hatX_x$, including those that induce embeddings $\bfC\into\hatcalO_{X,x}$ not landing in $\calO_{X,x}$.
\end{rem}

\subsection{Non-compact objects and hypersequences}
Because of functoriality of the algorithms from \S\ref{formal}, one immediately obtains a functorial
desingularization of non-compact qe schemes, formal schemes, and various analytic spaces. However, this time the
desingularization is just a projective morphism $X'\to X$ because its functorial splitting into a sequence of
blow ups can be infinite. One can often perform few blow ups simultaneously to obtain a finite splitting
$X'=X_n\longto X_0=X$, or an infinite splitting $$\dots\to X_n\to\dots\to X_1\to X_0=X$$ that reduces to a
finite sequence (and infinitely many empty blow ups) over any relatively compact subspace of $X$, but this is
not functorial. Moreover, despite some claims in the literature, even such an infinite splitting does not always
exist.

\begin{exam}
For concreteness, let us work complex-analytically and fix a strong desingularization algorithm $\calF$. It is
easy to construct a complex surface $S$ with an irreducible curve $C\subset S_\sing$ and a point $a_n\in S$ such
that the following conditions hold. At a generic point $\eta\in C$ one can locally describe $S$ by the equation
$f(x,y,z)=y^2+z^2=0$ (i.e. $S$ consists of two smooth branches meeting transversally along $C$), but $a_n$ is a
so special point that the resolution $\calF(S)$ on the appropriate charts over $a_n$ looks as follows: at least
$n$ times one blows up the preimage of $a_n$ on the strict transforms of $C$, and only then one blows up the
strict transform of $C$ (thus resolving the generic points of $C$). For example, an easy computation shows that
both for our algorithm $\calF$ and for the algorithm $\calF_\Var$ of Bierstone-Milman one can define the germ of
$S$ at $a_n$ by $y^2+z^2+x^{2n+3}=0$. Clearly, we can construct a non-compact surface with a curve $C$ and an
infinite sequence of points $a_1,a_2,\dots$ as above (use an infinite pasting procedure). If a factorization of
$\calF(S)$ into an infinite sequence $\dots\to S_1\to S_0=S$ would exist, then the strict transform of $C$ would
be a component of the center of some blow up, say $S_{n+1}\to S_n$. And this would contradict the assumption
that the composition is $\calF(S)$ over $a_{n+1}$.
\end{exam}

Nevertheless, there is a functorial way to split the desingularization. Instead of infinite blow up sequences
ordered by $\bfN$, one should consider their generalization, which is called blow up {\em hypersequences} in
\cite{temnon}. The latter are sequences ordered by a countable ordered set (the set of invariants of the
algorithm in this case) and such that over each relatively compact subobject $Y\into X$ the hypersequence
reduces to the finite blow up sequence $\calF(Y)$ saturated with infinitely many empty blow ups. Existence of
such a splitting is more or less a tautology and we refer to \cite[\S5.3]{temnon} for details.

\end{document}